\documentclass{amsart}
\usepackage{amssymb}
\usepackage{latexsym}
\usepackage{amsthm}
\theoremstyle{plain}

\newtheorem{thm}{Theorem}[section]

\newtheorem{ex}[thm]{Example}

\numberwithin{equation}{section}

\newcommand{\ZZ}{\mathbb Z}
\newcommand{\CC}{\mathbb C}

\newcommand{\RR}{\mathbb R}

\newcommand{\PP}{\mathbb P}
\newcommand{\fH}{\mathfrak H}

\begin{document}

\title{Explicit Equations of some Elliptic Modular 
surfaces} 
\author{Jaap Top}
\address{Department of Mathematics, University of Groningen,
P.O. Box 800, 9700 AV Groningen, The Netherlands}
\email{top@math.rug.nl}
\author{Noriko Yui}
\thanks{Noriko Yui was partially supported by a Research Grant
from NSERC, Canada.}
\address{Department of Mathematics and Statistics, Queen's University,
Kingston, Ontario Canada K7L 3N6}
\email{yui@mast.queensu.ca}
\date{\today}
\subjclass[2000]{14J27, 14J29}
\keywords{genus zero congruence subgroups,  
elliptic modular surfaces, surfaces of general type, cusp widths, 
semi-stable elliptic surface.}

\begin{abstract}
     We present explicit equations of semi-stable elliptic 
surfaces (i.e., having only type $I_n$ singular fibers)
which are associated to the torsion-free genus zero congruence 
subgroups of the modular group as classified by A. Sebbar.
\end{abstract}

\maketitle

\section{Introduction}
\label{sect1}

The purpose of this paper is to give explicit equations
for the elliptic modular surfaces associated
to torsion-free genus zero congruence subgroups $\Gamma$ of $PSL_2(\ZZ)$. 
We will call them {\it elliptic modular
surfaces of genus zero} for short. 
By the Noether formula, the Euler number of any elliptic surface
is a (nonnegative) multiple of $12$. It turns out that
an elliptic modular surface of genus zero has Euler number
one of $12, 24, 36, 48$ and $60$. Indeed, if it has Euler number
bigger than $12$, then it is semi-stable and the Euler number
equals the index of the corresponding group $\Gamma$
in $PSL_2(\ZZ)$. Sebbar's classification \cite{Se} of all such
torsion free congruence subgroups $\Gamma$ of genus zero implies that 
they have index $\leq 60$ in $PSL_2(\ZZ)$.

This paper can be regarded as a natural sequel to the article
 of Beauville \cite{Be}
which deals with elliptic modular surfaces (of genus zero) having
 Euler number $12$, and the
article of Livn\'e and Yui \cite{LY} which considered the
case of elliptic modular surfaces of genus zero with
Euler number $24$.

We will recall the results of Beauville \cite{Be} and of
Livn\'{e} and Yui \cite{LY} (presenting a different approach in the
latter case). Furthermore, we give explicit defining equations
for the elliptic modular surfaces with Euler number one of  
$36, 48$ and $60$.

Although our approach is different, these results could also be found
using McKay and Sebbar's tables in \cite{MS}, which provide one with the
$j$-invariant of each of the elliptic modular surfaces of genus zero.
We now illustrate how this works for a typical example (with Euler number $36$).

\begin{ex} Consider the group $\Gamma=\Gamma_0(2)\cap \Gamma(3)$. 
According to \cite[Table~5]{MS}, the
$j$-invariant $j(t)$ of the corresponding elliptic surface equals
\[
\frac{(t^3+4)^3(t^3+6t^2+4)^3(t^6-6t^5+36t^4+8t^3-24t^2+16)^3}
{t^6(t+1)^3(t^2-t+1)^3(t-2)^6(t^2+2t+4)^6}.
\]
To obtain a semi-stable elliptic surface with this $j$-invariant,
start from any elliptic surface with $j$-invariant $j$, such as
the one given by
\[
y^2+xy=x^3-\frac{36}{j-1728}x-\frac{1}{j-1728}.
\]
The corresponding elliptic surface over the $j$-line has
only three singular fibers, namely an $I_1$-fiber over $j=\infty$,
a $II$-fiber over $j=0$ and a $III^*$-fiber over $j=1728$.
Now substitute $j=j(t)$ in the  equation.
One obtains in this way an elliptic surface which
has $I_0^*$-fibers over all $t$ such that $j(t)=0$ (this can
be read off from the ramification of the function $j(t)$,
together with the fact that the surface given above has a fiber of type $II$
over $j=0$; compare \cite[p.~46]{BK}). Hence we have to take a quadratic twist over an
extension ramified at all these $t$ in order to make these fibers
smooth.
The equation $j(t)=\infty$ has four triple (one at $t=\infty$) and four six-fold roots,
hence the resulting surface has four fibers of type $I_3$ and four
fibers of type $I_6$. Finally, the equation $j(t)=1728$ 
has $18$ roots, each with multiplicity $2$. This implies that our
new surface has fibers of type $I_0^*$ over these $18$ values
of $t$. Therefore, after the quadratic twist over an extension
ramified only at these $18$ $t$'s, one obtains the  semi-stable
elliptic surface we were looking for.
\end{ex}

\medskip
It should be clear from the above example that this method
works in general; however, it is rather elaborate and the resulting
equations (even for a rather small example such as $j(t)$ above)
are quite cumbersome. Therefore, we offer an alternative
approach.

Our method consists of describing a different base change resulting
in the desired surface. Instead of starting from a surface with
$j$-invariant $j$, we start from a surface corresponding to
a torsion free group $\tilde{\Gamma}$ with 
$\Gamma\subset\tilde{\Gamma}\subset
PSL_2(\ZZ)$ and $\tilde{\Gamma}/\Gamma$ cyclic of order $>1$. The 
surprising
observation which makes our method work, is that such a group
exists for most cases in Sebbar's list.
We then start from the surface corresponding to $\tilde{\Gamma}$
and construct the one for $\Gamma$ as a cyclic base change of
the former.

We now illustrate this method on the same example as used above.

\begin{ex}\label{Gamma06}
The group $\Gamma:=\Gamma_0(2)\cap \Gamma(3)$ has index two in
$\tilde{\Gamma}:=\Gamma_0(6)$. The latter group already appears in the
paper of Beauville \cite{Be}; in fact, it corresponds to
the surface with affine equation
\[
(x+y)(y+1)(x+1)+txy=0.
\]
This surface has as singular fibers an $I_6$ (over $t=\infty$),
an $I_3$ (over $t=0$), an $I_2$ (over $t=1$)
and an $I_1$ (over $t=-8$). Now take the cyclic triple cover
of the $t$-line which ramifies over $t=1$ and $t=-8$ only.
Explicitly, this cover is described as
$s\mapsto t=(1+8s^3)/(1-s^3)$, for a coordinate $s$ on the cover
such that $s=0$ and $s=\infty$ are the branch points.
The resulting base changed surface has equation
\[
(s^3-1)(x+y)(y+1)(x+1)=(8s^3+1)xy,
\]
and its singular fibers are of type $I_6$ (over $s^3=1$ and
$s=0$) and of type $I_3$ (over $8s^3=-1$ and $s=\infty$).
\end{ex}

\medskip
We mention here that the problem of finding explicit
equations for elliptic surfaces (not necessarily modular)
has been getting considerable attention.  Miranda
and Persson \cite{MP} determined all possible configurations
of bad fibers of semi-stable elliptic K3 surfaces; for instance,
it turns out that such a surface must have at least $6$ (and at most
$24$) singular fibers, and there exist precisely $112$
possible configurations of bad fibers for such surfaces with
six singular fibers. These surfaces have Picard number $2+24-6=20$,
hence they are extremal elliptic K3 surfaces.
Such extremal surfaces have been studied by Nori \cite{No}.
Quite recently a short proof of a characterization of all such
surfaces was obtained by Kloosterman \cite{Kl03}.
W.~E.~Lang \cite{La1}, \cite{La2} studied rational extremal 
elliptic surfaces in positive characteristic.
The result of Miranda and Persson uses the construction of the
$j$-function of the surfaces involved. This relies on the Riemann
 existence theorem, hence it doesnot explicitly construct this
$j$-function. In particular, no equations of the surfaces are given.

The first semi-stable K3 surface in the list of Miranda and
Persson corresponds
to the configuration of singular fibers
$I_1,I_1,I_1,I_1,I_1,I_{19}$. Shioda \cite{Sh2} and Iron \cite{I} have
independently and with different methods obtained explicit equations 
realizing such a surface. In fact, over the complex numbers it is
unique: this follows from an observation of Kloosterman \cite[Thm.~8.2]{Kl03}
and also from the following argument due to Shioda \cite{Sh2}.
One may present any such surface by an equation
$y^2=x^3-3f(t)x-2g(t)$, with $f$ and $g$ polynomials of degree
$8$ resp. $12$ and `discriminant' $d(t):=g^2-f^3$ of degree $5$
(here the surface is given in such a way that the $I_{19}$-fiber
is over $t=\infty$). A special case of a result of Stothers \cite{St} shows that
upto obvious scalings, th equation $d=g^2-f^3$ has a unique
solution in polynomials of the given degrees. 

The paper of Livn\'{e} and Yui \cite{LY} presents equations of
nine other extremal semi-stable elliptic K3 surfaces, namely all
modular ones. In Section~\ref{2.3} below, it is explained how to
obtain such equations by the method of our paper. In fact, as we
will discuss in Section~\ref{remarks}, equations for at least $18$ of the $112$ entries
in the list of Miranda and Persson may be obtained from
Beauville's list mentioned above, by using quadratic base changes.

\section{Preliminaries}
\label{sect2}

\subsection{Congruence subgroups}
\label{cong}
Throughout, we use the notion `congruence subgroup' and the
notations $\Gamma(m)$, $\Gamma_1(m)$ and $\Gamma_0(m)$ in the
context of subgroups of $PSL_2(\ZZ)=SL_2(\ZZ)/\pm 1$. Note that
this differs slightly from a lot of literature where the same
notions and notations are used in the context of $SL_2(\ZZ)$;
we take the image under the natural quotient map of these.

Sebbar \cite{Se} gave the complete classification of
all torsion-free genus zero congruence subgroups of the modular
group $PSL_2(\ZZ)$.  
Since the groups are assumed to be torsion free, the natural morphism
from the corresponding
compact Riemann surface $\Gamma\backslash\fH^*$ (of genus zero) to
$PSL_2(\ZZ)\backslash \fH^*$ is ramified only in the cusps.
The ramification index is classically called the {\em cusp width}.
The Zeuthen-Hurwitz formula in this case gives
$\mu=6(k-2)$ where $\mu=[PSL_2{(\ZZ)}:\Gamma]$ is the index and
$k$ is the number of cusps for $\Gamma$. In particular, the index
is a multiple of $6$ for the type of subgroups considered.

In the table, the notation $\Gamma(m;m/d,\epsilon,\chi)$ is
used for a Larcher congruence subgroup. It is defined as follows.
Let $m>0$ be an integer and let $d>0$ be a divisor of $m$.
Write $m/d=h^2n$ with $n$ square-free. Take integers $\epsilon>0$ a divisor
of $h$, and $\chi>0$ a divisor of $\gcd(d\epsilon, m/(d\epsilon^2))$.
Then
\[
\Gamma(m;m/d,\epsilon,\chi):=\left\{\pm\left(\begin{array}{cc}1+\frac{m}{\epsilon\chi}\alpha
& d\beta\\ \frac{m}{\chi}\gamma&1+\frac{m}{\epsilon\chi}\delta\end{array}\right)
;\;\gamma\equiv\alpha\bmod\chi\right\}/\pm 1,
\]
where $\alpha,\beta,\gamma,\delta$ are integers and the matrices are supposed to
be in $SL_2(\ZZ)$.

The following table summarizes Sebbar's classification results. 
\medskip

\[
\begin{array}{|c|c|c|c|} \hline \hline
\mbox{Index} & t  & \mbox{Group} & \mbox{Cusp widths} \\ \hline
6            & 3  & \Gamma(2)    & 2,2,2              \\ 
             &    & \Gamma_0(4)  & 4,1,1              \\ \hline
12           & 4  & \Gamma(3)    & 3,3,3,3            \\
             &    & \Gamma_0(4)\cap \Gamma(2) & 4,4,2,2 \\
             &    & \Gamma_1(5)               & 5,5,1,1 \\
             &    & \Gamma_0(6)               & 6,3,2,1 \\
             &    & \Gamma_0(8)               & 8,2,1,1 \\
             &    & \Gamma_0(9)               & 9,1,1,1 \\ \hline
24           & 6  & \Gamma(4)                 & 4,4,4,4,4,4 \\
             &    & \Gamma_0(3)\cap\Gamma(2)  & 6,6,6,2,2,2 \\
             &    & \Gamma_1(7)               & 7,7,7,1,1,1 \\
             &    & \Gamma_1(8)               & 8,8,4,2,1,1 \\
             &    & \Gamma_0(8)\cap\Gamma(2)  & 8,8,2,2,2,2 \\
             &    & \Gamma(8;4,1,2)         & 8,4,4,4,2,2 \\
             &    & \Gamma_0(12)              & 12,4,3,3,1,1 \\
             &    & \Gamma_0(16)              & 16,4,1,1,1,1 \\
             &    & \Gamma(16;16,2,2)       & 16,2,2,2,1,1 \\ \hline
36           & 8  & \Gamma_0(2)\cap \Gamma(3) & 6,6,6,6,3,3,3,3 \\
             &    & \Gamma_1(9)               & 9,9,9,3,3,1,1,1 \\
             &    & \Gamma(9;3,1,3)         & 9,9,3,3,3,3,3,3 \\
             &    & \Gamma_1(10)              & 10,10,5,5,2,2,1,1 \\
             &    & \Gamma_0(18)              & 18,9,2,2,2,1,1,1 \\
             &    & \Gamma(27;27,3,3)       & 27,3,1,1,1,1,1,1 \\ \hline
48           & 10 & \Gamma_1(8)\cap\Gamma(2)  & 8,8,8,8,4,4,2,2,2,2 \\
             &    & \Gamma(8;2,1,2)         & 8,8,4,4,4,4,4,4,4,4 \\
             &    & \Gamma_1(12)              & 12,12,6,4,4,3,3,2,1,1 \\
             &    & \Gamma(12;6,1,2)        & 12,6,6,6,6,4,2,2,2,2 \\
             &    & \Gamma_0(16)\cap\Gamma_1(8) & 16,16,4,4,2,2,1,1,1,1 \\
             &    & \Gamma(16;8,2,2)        & 16,16,2,2,2,2,2,2,2,2 \\
             &    & \Gamma(24;24,2,2)       & 24,8,3,3,3,3,1,1,1,1 \\
             &    & \Gamma(32;32,4,2)       & 32,8,1,1,1,1,1,1,1,1 \\ \hline
60           & 12 & \Gamma(5)                 & 5,5,5,5,5,5,5,5,5,5,5,5 \\
             &    & \Gamma_0(25)\cap\Gamma_1(5) & 25,25,1,1,1,1,1,1,1,1,1,1\\
\hline \hline
\end{array}
\]
\centerline{Table~1: All torsion-free genus zero congruence subgroups}
\medskip







\subsection{Elliptic modular surfaces}\label{elmod}

Here we recall some basic facts about a special class of elliptic surfaces, called
{\it elliptic modular surfaces}, associated to congruence subgroups
$\Gamma\subset PSL_2(\ZZ)$ as discussed in Section~\ref{cong}.

The quotient of $\fH\times \CC$ by all automorphisms
of the form
$$(\tau, z)\mapsto (\gamma\tau, (c\tau +d)^{-1}(z+m\tau +n))$$
where $\gamma=\begin{pmatrix} a & b \\ c & d\end{pmatrix}\in\Gamma$
and $(m, n)\in\ZZ^2$ defines a surface equipped with a morphism to
the modular curve $X_\Gamma$ attached to $\Gamma$. The fiber over the image
in $X_\Gamma$ of a general
point $\tau\in \fH$  is the elliptic curve
corresponding to the lattice $\ZZ\oplus\ZZ\tau$. The surface 
obtained in this way can be extended to
an elliptic surface over the modular curve $X_\Gamma$.  This surface is
called the {\it elliptic modular surface} associated to $\Gamma$.

McKay and Sebbar \cite{MS} observed
that all groups $\Gamma$ in the list given in subsection~\ref{cong}
(except the two which have index $6$ in $PSL_2(\ZZ)$)
have an explicit lift $\overline{\Gamma}$ to $SL_2(\ZZ)$  with the following properties:
\begin{itemize}
\item[(a)] $\overline{\Gamma}$ has no elliptic elements,
\item[(b)]  $\overline{\Gamma}$ contains no element of trace equal to $-2$.  
\end{itemize}
Hence following  a classical method of Kodaira, they obtained
the following result.
\begin{thm}\label{thm2.1} {\sl Let $\Gamma$ be one of the subgroups in Table 1
with index $\mu>6$ and let $t$ be the number of inequivalent cusps for $\Gamma$.  Then
the elliptic modular surfaces attached to $\Gamma$ is
semi-stable, and has $t$ singular fibers of type $I_n$ corresponding
to these cusps. Moreover, the fiber over such a cusp is of type $I_n$
where $n$ equals the cusp width.  }
\end{thm}

Moreover, McKay and Sebbar explicitly give the relation between the Hauptmodul $t$ for $\Gamma$ (which
is by definition a generator of the function field of the rational
curve $X_\Gamma$) and the
elliptic modular function $j$. In other words, they present the rational function
$j(t)$ which gives the $j$-invariant of the corresponding elliptic modular surface.
Let $I_{n_1},\cdots, I_{n_k}$ be the singular fibers of such
a surface, then Kodaira's results imply
$\sum n_i=\mu$.
In particular,
the index $\mu$ gives the Euler number of the surface, which in turn
gives its geometric genus $p_g=\mu/12-1$.  
Therefore, the index $12$ subgroups
give rise to rational elliptic surfaces, the index $24$
to elliptic K3 surfaces, and the remaining indices $36, 48$
and $60$ to elliptic surfaces of general type with $p_g=2,3$ and
$4$, respectively. 

These elliptic surfaces are {\em extremal}, which means that their N\'{e}ron-Severi
group has rank equal to their Hodge number $h^{1,1}$ and is generated by
a section and all components of singular fibers. To see this, note
that the Euler number equals the sum of Hodge numbers
$h^0+h^{0,2}+h^{1,1}+h^{2,0}+h^4=2+2p_g+h^{1,1}$,
hence $h^{1,1}=5\mu/6$. The components of singular fibers plus a section
generate an indivisible subgroup of the  N\'{e}ron-Severi group of rank
$2+\sum_{i=1}^k(n_i-1)=2+\mu-k$. In Section~\ref{cong} it was remarked that
$\mu=6(k-2)$. Therefore, it follows that $2+\mu-k=5\mu/6=h^{1,1}$,
which shows that these surfaces are extremal.

\subsection{Defining equations corresponding to index
$12$ (resp. $24$)}\label{2.3}

The defining equations for the rational elliptic surfaces
were given by Beauville \cite{Be} (not in a Weierstrass form; however,
it is easy to derive Weierstrass equations from the given ones as is explained,
e.g., in \cite[Chapter~8]{Ca}).
Equations for the elliptic K3 surfaces were obtained by Livn\'e and Yui \cite{LY}
(in a Weierstrass form).

For the sake of completeness we list these known cases. First
we list Beauville's defining equations
for the rational elliptic modular surfaces in $\PP^2\times \PP^1$
where $[x:y:z]$ denotes the projective coordinate for $\PP^2$
and $t\in\PP^1$ is the affine coordinate for $\PP^1$.
\medskip

\[
\begin{array}{|c|c|c|} \hline \hline
\mbox{modular curve}&\mbox{group} & \mbox{defining equation} \\ \hline 
X(3)  &\Gamma(3) & x^3+y^3+z^3+t\,xyz=0  \\ 
      &\Gamma_0(4)\cap \Gamma(2) & x(x^2+z^2+2zy)+t\,z(x^2-y^2)=0 \\ 
X_1(5)&\Gamma_1(5) & x(x-z)(y-z)+t\,zy(x-y)=0  \\  
X_0(6)&\Gamma_0(6) & (x+y)(y+z)(z+x)+t\,xyz=0  \\  
X_0(8)&\Gamma_0(8) & (x+y)(xy-z^2)+t\,xyz=0   \\  
X_0(9)&\Gamma_0(9) & x^2y+y^2z+z^2x+t\, xyz=0 \\ \hline \hline 
\end{array}
\]
\centerline{Table 2: Defining equations corresponding to index $12$ groups}
\medskip

Now we list defining equations for the elliptic modular
K3 surfaces obtained in Livn\'e and Yui \cite{LY}. As 
Fuchsian groups in $PSL_2(\RR)$, the nine groups 
of index $24$ are further classified into four sets of
mutually conjugate groups:

$$\Gamma(4)\;\mbox{ and } \Gamma_0(8)\cap \Gamma(2)\;\mbox{ and }\Gamma_0(16),$$ 
$$\Gamma_0(3)\cap \Gamma(2)\;\mbox{ and } \Gamma_0(12),$$
$$\Gamma_1(7),$$
$$\Gamma_1(8),\;\mbox{ and } \Gamma_1(8;4,1,2)\;\mbox{ and } \Gamma_1(16;16,2,2).$$ 

Conjugation of these groups translates into isogenies between the
corresponding elliptic surfaces. In fact, all isogenies which show up here,
are composed of isogenies of degree $2$. Since formulas for the image under such
isogenies are well known \cite[pp.~58-59]{Ca}, we only present one equation
from each isogeny class.
\medskip

\[
\begin{array}{|c|c|} \hline \hline
\mbox{group} & \mbox{defining equation} \\ \hline
\Gamma(4)        & y^2=x(x-1)(x-(t+t^{-1})^2/4)
                   \\ \hline 
\Gamma_1(7) & y^2+(1+t-t^2)xy+(t^2-t^3)y 
             =x^3+(t^2-t^3)x  \\ \hline
\Gamma(8;4,1,2) & y^2=x^3-2(8t^4-16t^3+16t^2-8t+1)x^2
                +(8t^2-8t+1)(2t-1)^4x   \\ \hline 
\Gamma_0(12) & y^2+(t^2+1)xy-t^2(t^2-1)y 
             =x^3-t^2(t^2-1)x^2  \\ \hline
\hline
\end{array}
\]
\centerline{Table 3: Defining equations corresponding to index $24$ groups}
\medskip

We now describe a different method to obtain defining
equations corresponding to all nine index $24$ groups. 
In fact, this geometric method works not only for index $24$
groups but for other groups as well. 
This will be done by realizing all but one of these elliptic modular surfaces
as double covers of  Beauville's rational elliptic modular 
surfaces and then by deriving defining equations from those of
Beauville. The remaining $9$th surface is the one corresponding
to $\Gamma_1(7)$. An explicit equation for this case was already
obtained by Kubert \cite{Ku} in 1976; see also \cite[Table~3]{HLP}.
It is given by
\[
y^2+(1-c)xy-by=x^3-bx^2
\]
where $b=t^3-t^2$ and $c=t^2-t$.

\subsubsection{Two subgroups of $\Gamma_0(4)\cap\Gamma(2)$}\label{24-1}\mbox{ }\\
The first example we consider, is the index $12$ subgroup
$\Gamma_0(4)\cap\Gamma(2)$. According to Beauville, it
corresponds to the surface with affine equation
\[
x(x^2+2y+1)+t(x^2-y^2)=0.
\]
A Weierstrass equation for the same surface is obtained by
taking $t$ and $\xi:=tx$ and $\eta:=t^2y-tx$ as new coordinates. The
resulting equation is
\[ \eta^2=\xi^3+(t^2+1)\xi^2+t^2\xi.
\]
The configuration of singular fibers of this surface is
given in the following table.

\medskip
\begin{center}
\begin{tabular}{|r|cccc|}\hline\hline
$t=$ & $0$ & $1$ & $-1$ & $\infty$ \\\hline
fiber: & $I_4$ & $I_2$ & $I_2$ & $I_4$\\ \hline\hline
\end{tabular}
\end{center}

\medskip
The group $\Gamma_0(4)\cap\Gamma(2)$ contains $\Gamma(4)$ as well as
$\Gamma(8;4,1,2)$ as index $2$ subgroups.
The first one corresponds to a double cover of modular curves
ramified over the cusps of width $2$. The second one yields a
double cover ramified over one cusp of width $2$ and one cusp of
width $4$.

An explicit rational map which realizes the first double cover, is
\[s\mapsto t:=(s^2-1)/(s^2+1),
\]
with $s=0$ and $s=\infty$ as ramification points.
Hence the base changed surface with equation
\[
(s^2+1)x(x^2+2y+1)+(s^2-1)(x^2-y^2)=0
\]
is a semi-stable surface with $6$ singular fibers of type $I_4$,
corresponding to $\Gamma(4)$. We express the data of the
associated double cover of modular curves by the following diagram.

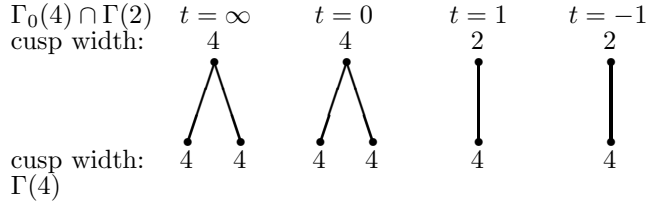
\vspace{\baselineskip}

\begin{center}
\begin{picture}(230,70)(-30,0)
\thicklines
\put(-27,65){$\Gamma_0(4)\cap\Gamma(2)$}
\put(-27,55){cusp width:}

\put(147,55){$2$} 
\put(139,65){$t=1$}
\put(150,50){\circle*{3}}
\put(150,50){\line(0,-1){30}}
\put(150,20){\circle*{3}}
\put(147,10){$4$}

\put(197,55){$2$}
\put(184,65){$t=-1$}
\put(200,50){\circle*{3}}
\put(200,50){\line(0,-1){30}}
\put(200,20){\circle*{3}}
\put(197,10){$4$}

\put(47,55){$4$}
\put(37,65){$t=\infty$}
\put(50,50){\circle*{3}}
\put(50,50){\line(-1,-3){10}}
\put(50,50){\line(1,-3){10}}
\put(40,20){\circle*{3}}
\put(60,20){\circle*{3}}
\put(37,10){$4$}
\put(57,10){$4$}

\put(97,55){$4$}
\put(88,65){$t=0$}
\put(100,50){\circle*{3}}
\put(100,50){\line(-1,-3){10}}
\put(100,50){\line(1,-3){10}}
\put(90,20){\circle*{3}}
\put(110,20){\circle*{3}}
\put(87,10){$4$}
\put(107,10){$4$}

\put(-27,10){cusp width:}
\put(-27,0){$\Gamma(4)$}
\end{picture}
\end{center}

\vspace{\baselineskip}

Similarly, an explicit rational map realizing the second double cover
is 
\[
s\mapsto t:= s^2+1.
\]

This corresponds to the diagram

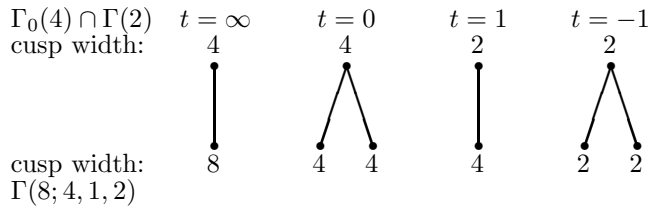
\vspace{\baselineskip}

\begin{center}
\begin{picture}(230,70)(-30,0)
\thicklines
\put(-27,65){$\Gamma_0(4)\cap\Gamma(2)$}
\put(-27,55){cusp width:}

\put(47,55){$4$} 
\put(37,65){$t=\infty$}
\put(50,50){\circle*{3}}
\put(50,50){\line(0,-1){30}}
\put(50,20){\circle*{3}}
\put(47,10){$8$}

\put(147,55){$2$}
\put(139,65){$t=1$}
\put(150,50){\circle*{3}}
\put(150,50){\line(0,-1){30}}
\put(150,20){\circle*{3}}
\put(147,10){$4$}

\put(197,55){$2$}
\put(185,65){$t=-1$}
\put(200,50){\circle*{3}}
\put(200,50){\line(-1,-3){10}}
\put(200,50){\line(1,-3){10}}
\put(190,20){\circle*{3}}
\put(210,20){\circle*{3}}
\put(187,10){$2$}
\put(207,10){$2$}

\put(97,55){$4$}
\put(89,65){$t=0$}
\put(100,50){\circle*{3}}
\put(100,50){\line(-1,-3){10}}
\put(100,50){\line(1,-3){10}}
\put(90,20){\circle*{3}}
\put(110,20){\circle*{3}}
\put(87,10){$4$}
\put(107,10){$4$}

\put(-27,10){cusp width:}
\put(-27,0){$\Gamma(8;4,1,2)$}
\end{picture}
\end{center}

\vspace{\baselineskip}

Hence $\Gamma(8;4,1,2)$ is realized by
\[
x(x^2+2y+1)+(s^2+1)(x^2-y^2)=0.
\]

\subsubsection{Two subgroups of $\Gamma_0(6)$}\label{24-2}\mbox{ }\\
The Beauville surface corresponding to $\Gamma_0(6)$ was
already described in Example~\ref{Gamma06}.
Its equation is 
\[(x+y)(x+1)(y+1)+txy=0
\] and its
singular fiber configuration is given by

\medskip
\begin{center}
\begin{tabular}{|r|cccc|}\hline\hline
$t=$ & $-8$ & $1$ & $0$ & $\infty$ \\\hline
fiber: & $I_1$ & $I_2$ & $I_3$ & $I_6$\\ \hline\hline
\end{tabular}
\end{center}

\medskip
Two subgroups of index $2$ of $\Gamma_0(6)$ are $\Gamma_0(3)\cap \Gamma(2)$
and $\Gamma_0(12)$. The corresponding double cover of modular curves
ramifies over the cusps of width $1$ and $3$ in the first case.
In the second case, it ramifies over the cusps of width $2$ and $6$.

So to realize the surface corresponding to $\Gamma_0(3)\cap \Gamma(2)$
one can use the base change obtained from
\[
s\mapsto t:=s^2/(1-8s^2).
\]
This yields the diagram

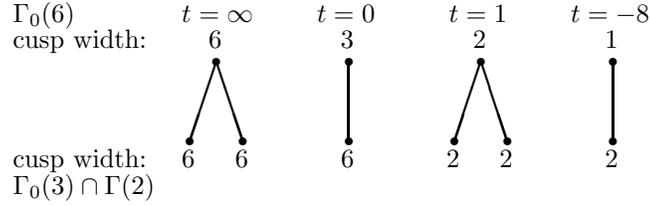
\vspace{\baselineskip}

\begin{center}
\begin{picture}(230,70)(-30,0)
\thicklines
\put(-27,65){$\Gamma_0(6)$}
\put(-27,55){cusp width:}

\put(97,55){$3$} 
\put(88,65){$t=0$}
\put(100,50){\circle*{3}}
\put(100,50){\line(0,-1){30}}
\put(100,20){\circle*{3}}
\put(97,10){$6$}

\put(197,55){$1$}
\put(184,65){$t=-8$}
\put(200,50){\circle*{3}}
\put(200,50){\line(0,-1){30}}
\put(200,20){\circle*{3}}
\put(197,10){$2$}

\put(47,55){$6$}
\put(37,65){$t=\infty$}
\put(50,50){\circle*{3}}
\put(50,50){\line(-1,-3){10}}
\put(50,50){\line(1,-3){10}}
\put(40,20){\circle*{3}}
\put(60,20){\circle*{3}}
\put(37,10){$6$}
\put(57,10){$6$}

\put(147,55){$2$}
\put(138,65){$t=1$}
\put(150,50){\circle*{3}}
\put(150,50){\line(-1,-3){10}}
\put(150,50){\line(1,-3){10}}
\put(140,20){\circle*{3}}
\put(160,20){\circle*{3}}
\put(137,10){$2$}
\put(157,10){$2$}

\put(-27,10){cusp width:}
\put(-27,0){$\Gamma_0(3)\cap\Gamma(2)$}
\end{picture}
\end{center}

\vspace{\baselineskip}

For the group $\Gamma_0(12)$ the map
\[
s\mapsto t:=1-s^2
\]
gives the desired covering. In a diagram:

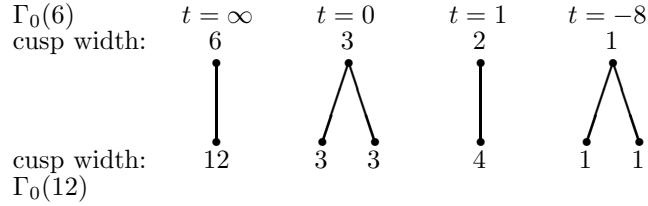
\vspace{\baselineskip}

\begin{center}
\begin{picture}(230,70)(-30,0)
\thicklines
\put(-27,65){$\Gamma_0(6)$}
\put(-27,55){cusp width:}

\put(47,55){$6$} 
\put(37,65){$t=\infty$}
\put(50,50){\circle*{3}}
\put(50,50){\line(0,-1){30}}
\put(50,20){\circle*{3}}
\put(45,10){$12$}

\put(147,55){$2$}
\put(138,65){$t=1$}
\put(150,50){\circle*{3}}
\put(150,50){\line(0,-1){30}}
\put(150,20){\circle*{3}}
\put(147,10){$4$}

\put(97,55){$3$}
\put(88,65){$t=0$}
\put(100,50){\circle*{3}}
\put(100,50){\line(-1,-3){10}}
\put(100,50){\line(1,-3){10}}
\put(90,20){\circle*{3}}
\put(110,20){\circle*{3}}
\put(87,10){$3$}
\put(107,10){$3$}

\put(197,55){$1$}
\put(183,65){$t=-8$}
\put(200,50){\circle*{3}}
\put(200,50){\line(-1,-3){10}}
\put(200,50){\line(1,-3){10}}
\put(190,20){\circle*{3}}
\put(210,20){\circle*{3}}
\put(187,10){$1$}
\put(207,10){$1$}

\put(-27,10){cusp width:}
\put(-27,0){$\Gamma_0(12)$}
\end{picture}
\end{center}

\vspace{\baselineskip}

Hence one obtains the surfaces
\[
\mbox{ for }\;\Gamma_0(3)\cap \Gamma(2):\quad 
(8s^2-1)(x+y)(x+1)(y+1)=s^2xy
\]
and
\[
\mbox{ for }\;\Gamma_0(12):\quad
(x+y)(x+1)(y+1)+(1-s^2)xy=0.
\]

\subsubsection{Four subgroups of $\Gamma_0(8)$}\label{24-3}\mbox{ }\\
The Beauville surface corresponding to $\Gamma_0(8)$
was given by Beauville as 
\[
(x+y)(xy-1)+\tilde{t}xy=0.
\] 
Alternatively, it is defined
by the Weierstrass equation
\[
\eta^2=\xi^3+(2-{t}^2)\xi^2+\xi,
\] 
in which
$\eta=(x+y-\frac12txy)y^{-2}$ and $\xi=x/y$ and $t=i\tilde{t}/2$. 
The configuration
of singular fibers of this surface is as follows.

\medskip
\begin{center}
\begin{tabular}{|r|cccc|}\hline\hline
$t=$ & $\infty$ & $0$ & $1$ & $-1$ \\\hline
fiber: & $I_8$ & $I_2$ & $I_1$ & $I_1$\\ \hline\hline
\end{tabular}
\end{center}

\medskip
The reason for rescaling Beauville's parameter written as $\tilde{t}$ here,
is precisely the fact that in the new parameter $t$, all singular fibers
occur over rational points.
The group $\Gamma_0(8)$ turns out to contain four index $2$
genus zero congruence subgroups, namely
\[
\Gamma_1(8),\quad \Gamma_0(8)\cap\Gamma(2),\quad \Gamma_0(16),\quad\mbox{and}
\quad
\Gamma(16;16,2,2).
\]
The corresponding double covers of modular curves ramify over
two cusps, which are of width $1$ and $2$ in the first case,
both of width $1$ in the second case, of width $8$ and $2$
in the third case, and of width $8$ and $1$ in the remaining case.
This is expressed in the next four diagrams.

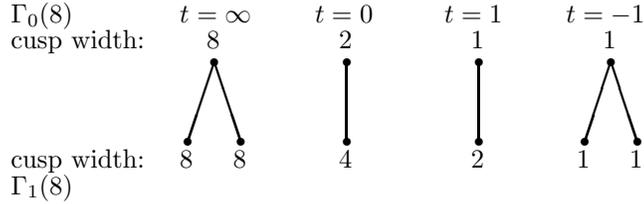
\vspace{\baselineskip}

\begin{center}
\begin{picture}(230,70)(-30,0)
\thicklines
\put(-27,65){$\Gamma_0(8)$}
\put(-27,55){cusp width:}

\put(147,55){$1$} 
\put(137,65){$t=1$}
\put(150,50){\circle*{3}}
\put(150,50){\line(0,-1){30}}
\put(150,20){\circle*{3}}
\put(147,10){$2$}

\put(97,55){$2$}
\put(88,65){$t=0$}
\put(100,50){\circle*{3}}
\put(100,50){\line(0,-1){30}}
\put(100,20){\circle*{3}}
\put(97,10){$4$}

\put(47,55){$8$}
\put(37,65){$t=\infty$}
\put(50,50){\circle*{3}}
\put(50,50){\line(-1,-3){10}}
\put(50,50){\line(1,-3){10}}
\put(40,20){\circle*{3}}
\put(60,20){\circle*{3}}
\put(37,10){$8$}
\put(57,10){$8$}

\put(197,55){$1$}
\put(183,65){$t=-1$}
\put(200,50){\circle*{3}}
\put(200,50){\line(-1,-3){10}}
\put(200,50){\line(1,-3){10}}
\put(190,20){\circle*{3}}
\put(210,20){\circle*{3}}
\put(187,10){$1$}
\put(207,10){$1$}

\put(-27,10){cusp width:}
\put(-27,0){$\Gamma_1(8)$}
\end{picture}
\end{center}

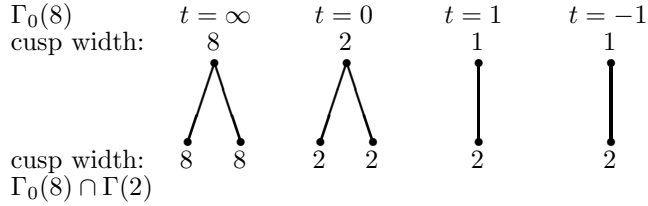
\vspace{\baselineskip}

\begin{center}
\begin{picture}(230,70)(-30,0)
\thicklines
\put(-27,65){$\Gamma_0(8)$}
\put(-27,55){cusp width:}

\put(147,55){$1$} 
\put(137,65){$t=1$}
\put(150,50){\circle*{3}}
\put(150,50){\line(0,-1){30}}
\put(150,20){\circle*{3}}
\put(147,10){$2$}

\put(197,55){$1$}
\put(184,65){$t=-1$}
\put(200,50){\circle*{3}}
\put(200,50){\line(0,-1){30}}
\put(200,20){\circle*{3}}
\put(197,10){$2$}

\put(47,55){$8$}
\put(37,65){$t=\infty$}
\put(50,50){\circle*{3}}
\put(50,50){\line(-1,-3){10}}
\put(50,50){\line(1,-3){10}}
\put(40,20){\circle*{3}}
\put(60,20){\circle*{3}}
\put(37,10){$8$}
\put(57,10){$8$}

\put(97,55){$2$}
\put(88,65){$t=0$}
\put(100,50){\circle*{3}}
\put(100,50){\line(-1,-3){10}}
\put(100,50){\line(1,-3){10}}
\put(90,20){\circle*{3}}
\put(110,20){\circle*{3}}
\put(87,10){$2$}
\put(107,10){$2$}

\put(-27,10){cusp width:}
\put(-27,0){$\Gamma_0(8)\cap\Gamma(2)$}
\end{picture}
\end{center}

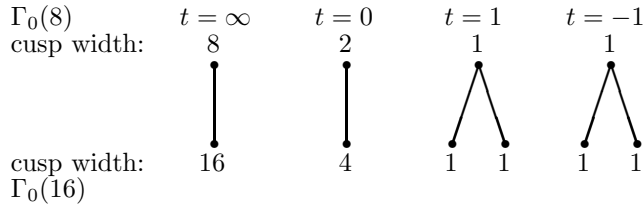
\vspace{\baselineskip}

\begin{center}
\begin{picture}(230,70)(-30,0)
\thicklines
\put(-27,65){$\Gamma_0(8)$}
\put(-27,55){cusp width:}

\put(47,55){$8$} 
\put(37,65){$t=\infty$}
\put(50,50){\circle*{3}}
\put(50,50){\line(0,-1){30}}
\put(50,20){\circle*{3}}
\put(44,10){$16$}

\put(97,55){$2$}
\put(88,65){$t=0$}
\put(100,50){\circle*{3}}
\put(100,50){\line(0,-1){30}}
\put(100,20){\circle*{3}}
\put(97,10){$4$}

\put(147,55){$1$}
\put(137,65){$t=1$}
\put(150,50){\circle*{3}}
\put(150,50){\line(-1,-3){10}}
\put(150,50){\line(1,-3){10}}
\put(140,20){\circle*{3}}
\put(160,20){\circle*{3}}
\put(137,10){$1$}
\put(157,10){$1$}

\put(197,55){$1$}
\put(183,65){$t=-1$}
\put(200,50){\circle*{3}}
\put(200,50){\line(-1,-3){10}}
\put(200,50){\line(1,-3){10}}
\put(190,20){\circle*{3}}
\put(210,20){\circle*{3}}
\put(187,10){$1$}
\put(207,10){$1$}

\put(-27,10){cusp width:}
\put(-27,0){$\Gamma_0(16)$}
\end{picture}
\end{center}

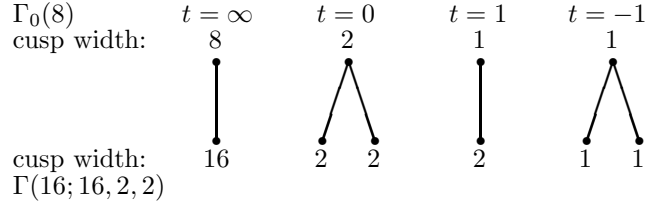
\vspace{\baselineskip}

\begin{center}
\begin{picture}(230,70)(-30,0)
\thicklines
\put(-27,65){$\Gamma_0(8)$}
\put(-27,55){cusp width:}

\put(47,55){$8$} 
\put(37,65){$t=\infty$}
\put(50,50){\circle*{3}}
\put(50,50){\line(0,-1){30}}
\put(50,20){\circle*{3}}
\put(45,10){$16$}

\put(147,55){$1$}
\put(138,65){$t=1$}
\put(150,50){\circle*{3}}
\put(150,50){\line(0,-1){30}}
\put(150,20){\circle*{3}}
\put(147,10){$2$}

\put(97,55){$2$}
\put(88,65){$t=0$}
\put(100,50){\circle*{3}}
\put(100,50){\line(-1,-3){10}}
\put(100,50){\line(1,-3){10}}
\put(90,20){\circle*{3}}
\put(110,20){\circle*{3}}
\put(87,10){$2$}
\put(107,10){$2$}

\put(197,55){$1$}
\put(183,65){$t=-1$}
\put(200,50){\circle*{3}}
\put(200,50){\line(-1,-3){10}}
\put(200,50){\line(1,-3){10}}
\put(190,20){\circle*{3}}
\put(210,20){\circle*{3}}
\put(187,10){$1$}
\put(207,10){$1$}

\put(-27,10){cusp width:}
\put(-27,0){$\Gamma(16;16,2,2)$}
\end{picture}
\end{center}

\vspace{\baselineskip}

Explicit covering maps which can be used for this, are presented
in the following table.
\[
\begin{array}{rcl}
\mbox{ for }\;\Gamma_1(8): & & s\mapsto t:=s^2/(s^2+1);\\
\mbox{ for }\;\Gamma_0(8)\cap\Gamma(2): & & s\mapsto t:=(s^2-1)/(s^2+1);\\
\mbox{ for }\;\Gamma_0(16): & & s\mapsto t:=s^2;\\
\mbox{ for }\;\Gamma(16;16,2,2): & &s\mapsto t:=s^2+1.
\end{array}
\]
It should be evident how to obtain explicit equations for the
base changed surfaces from this.

\section{The main result}
\label{sect3}

In this section we will give explicit equations
for the elliptic surfaces of general type corresponding
to the remaining genus zero congruence subgroups of
index $36, 48$ and $60$ in
Table 1.
The observation which makes this a rather easy task,
is the following.
\begin{thm}\label{main}
With the exception of the two groups $\Gamma_1(7)$ and $\Gamma_1(10)$,
all torsion-free genus zero congruence subgroups of index at
least $24$ in $PSL_2(\ZZ)$ correspond to modular curves which
are obtained by a composition of cyclic covers from the
modular curves in Beauville's list (see Table~2).

Here, a composition of cyclic covers is only needed for the groups
$\Gamma_1(12)$ and $\Gamma_0(16)\cap \Gamma_1(8)$; in all other
cases, a single cyclic map suffices.
\end{thm}

\medskip
For the case of index $24$, Theorem~\ref{main} was proven
in the Paragraphs~\ref{24-1}, \ref{24-2}, and \ref{24-3} above.
We will consider the remaining cases in Section~\ref{proof} below.

Note that the exceptional case $\Gamma_1(7)$ was already discussed
in Section~\ref{2.3}. For the remaining exceptional case $\Gamma_1(10)$,
we quote \cite[Table~6]{HLP}, which provides the equation
\[y^2=x(x^2+ax+b)
\]
with
\[
a=-(2t^2-2t+1)(4t^4-12t^3+6t^2+2t-1)
\]
and
\[
b=16(t^2-3t+1)(t-1)^5t^5
\]
for the corresponding elliptic modular surface.
This is a surface with two $I_5$-fibers at $t=\infty$ and at $t=1/2$,
two $I_{10}$-fibers at $t=0$ and at $t=1$, two $I_2$-fibers at the
roots of $t^2-3t+1=0$, and finally, two $I_1$-fibers at the roots of
$4t^2-2t-1=0$.

\medskip
Note that a cyclic cover $\varphi:\PP^1\to\PP^1$
of degree $d$ ramifies over exactly two points. Write these ramification
points as $t=a$ and $t=b$ for some coordinate $t$ on the
second $\PP^1$ (which, for convenience, we take in such a way
that $a\neq\infty$). Choose a coordinate $s$ on the first $\PP^1$
 such that
$s=0$ is mapped to $t=a$ and $s=\infty$ is mapped to $t=b$.
With such choices, $\varphi$ can be given as
\[
\varphi:\;\; s\mapsto t:=\frac{bs^d+\lambda a}{s^d+\lambda}
\]
for an arbitrary choice of $\lambda\neq 0$. In case $b=\infty$
one can take $\varphi(s)=\lambda s^d+a$.
Hence once we have determined the degrees and the ramification points
for the maps needed in Theorem~\ref{main}, it is evident how
to obtain explicit equations for the corresponding elliptic
modular surfaces starting from Beauville's equations.

\section{Proof of the main result}\label{proof}
Here we present a proof of Theorem~\ref{main}, by considering the
necessary groups in Beauville's list one by one.
\subsection{An index $3$  subgroup of $\Gamma(3)$}\mbox{ }\\
The elliptic modular surface associated with $\Gamma(3)$ is the
famous Hessian family given in affine form as
\[
x^3+y^3+1-3txy=0.
\]
Note that we have slightly rescaled the parameter used by Beauville.
The family has four fibers of type $I_3$, corresponding to
$t=\infty,1,\omega,\bar{\omega}$. Here $\omega$ denotes a primitive
cube root of unity (so $\omega^2+\omega+1=0$).
The subgroup $\Gamma(9;3,1,3)$ yields a cyclic covering of degree $3$
of the modular curve $X(3)$. It ramifies over two of the cusps.
This is illustrated in the following diagram.

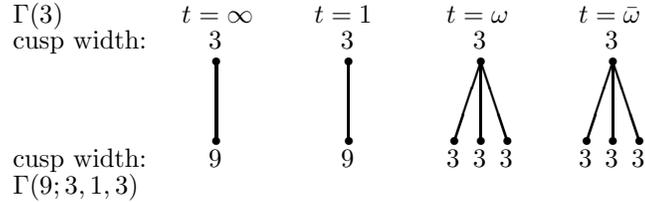
\vspace{\baselineskip}

\begin{center}
\begin{picture}(230,70)(-30,0)
\thicklines
\put(-27,65){$\Gamma(3)$}
\put(-27,55){cusp width:}

\put(47,55){$3$} 
\put(37,65){$t=\infty$}
\put(50,50){\circle*{3}}
\put(50,50){\line(0,-1){30}}
\put(50,20){\circle*{3}}
\put(47,10){$9$}

\put(97,55){$3$}
\put(87,65){$t=1$}
\put(100,50){\circle*{3}}
\put(100,50){\line(0,-1){30}}
\put(100,20){\circle*{3}}
\put(97,10){$9$}

\put(147,55){$3$}
\put(137,65){$t=\omega$}
\put(150,50){\circle*{3}}
\put(150,50){\line(-1,-3){10}}
\put(150,50){\line(0,-1){30}}
\put(150,50){\line(1,-3){10}}
\put(140,20){\circle*{3}}
\put(150,20){\circle*{3}}
\put(160,20){\circle*{3}}
\put(137,10){$3$}
\put(147,10){$3$}
\put(157,10){$3$}

\put(197,55){$3$}
\put(187,65){$t=\bar{\omega}$}
\put(200,50){\circle*{3}}
\put(200,50){\line(-1,-3){10}}
\put(200,50){\line(1,-3){10}}
\put(200,50){\line(0,-1){30}}
\put(190,20){\circle*{3}}
\put(200,20){\circle*{3}}
\put(210,20){\circle*{3}}
\put(187,10){$3$}
\put(197,10){$3$}
\put(207,10){$3$}

\put(-27,10){cusp width:}
\put(-27,0){$\Gamma(9;3,1,3)$}
\end{picture}
\end{center}

\subsection{Two index $3$ subgroups of $\Gamma_0(9)$}\mbox{ }\\
The elliptic modular surface associated with $\Gamma_0(9)$ is
presented by Beauville as
\[
x^2y+y^2+x+txy=0.
\]
Using new coordinates $\xi:=-4/x$ and $\eta:=(8y+4x^2+4tx)/x^2$
this is transformed into the Weierstrass equation
\[
\eta^2=\xi^3+t^2\xi^2-8t\xi+16.
\]
This surface has three singular fibers of type $I_1$
at $t=-3,-3\omega,-3\bar{\omega}$ (with $\omega^2+\omega+1=0$)
and a fiber of type $I_9$ at $t=\infty$.

The subgroup $\Gamma_1(9)$ corresponds to a cyclic triple covering
of $X_0(9)$ ramified over two cusps of width $1$.
With the parameter $t$ used here, this covering can be expressed by
\[
s\mapsto t:=\frac{s^3-3s+1}{s^2-s}.
\]
It corresponds to the following diagram.

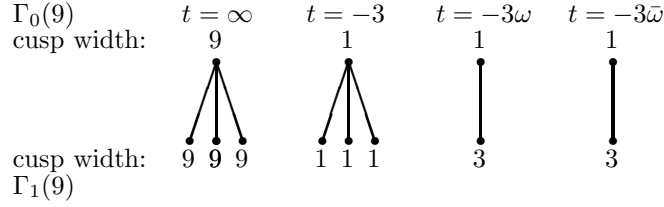
\vspace{\baselineskip}

\begin{center}
\begin{picture}(230,70)(-30,0)
\thicklines
\put(-27,65){$\Gamma_0(9)$}
\put(-27,55){cusp width:}

\put(47,55){$9$} 
\put(37,65){$t=\infty$}
\put(50,50){\circle*{3}}
\put(50,50){\line(0,-1){30}}
\put(50,20){\circle*{3}}
\put(47,10){$9$}
\put(50,50){\line(-1,-3){10}}
\put(50,50){\line(0,-1){30}}
\put(50,50){\line(1,-3){10}}
\put(40,20){\circle*{3}}
\put(50,20){\circle*{3}}
\put(60,20){\circle*{3}}
\put(37,10){$9$}
\put(47,10){$9$}
\put(57,10){$9$}

\put(97,55){$1$}
\put(84,65){$t=-3$}
\put(100,50){\circle*{3}}
\put(100,50){\line(-1,-3){10}}
\put(100,50){\line(0,-1){30}}
\put(100,50){\line(1,-3){10}}
\put(90,20){\circle*{3}}
\put(100,20){\circle*{3}}
\put(110,20){\circle*{3}}
\put(87,10){$1$}
\put(97,10){$1$}
\put(107,10){$1$}

\put(147,55){$1$}
\put(133,65){$t=-3\omega$}
\put(150,50){\circle*{3}}
\put(150,50){\line(0,-1){30}}
\put(150,20){\circle*{3}}
\put(147,10){$3$}

\put(197,55){$1$}
\put(183,65){$t=-3\bar{\omega}$}
\put(200,50){\circle*{3}}
\put(200,50){\line(0,-1){30}}
\put(200,20){\circle*{3}}
\put(197,10){$3$}

\put(-27,10){cusp width:}
\put(-27,0){$\Gamma_1(9)$}
\end{picture}
\end{center}

\vspace{\baselineskip}

The subgroup $\Gamma(27;27,3,3)$ yields the cycclic triple covering of
$X_0(9)$ ramifying over the other two cusps (of width $9$ and $1$,
respectively). Hence in a diagram:

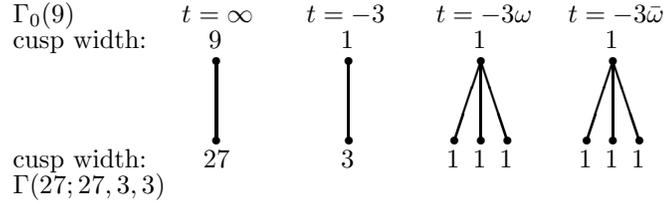
\vspace{\baselineskip}

\begin{center}
\begin{picture}(230,70)(-30,0)
\thicklines
\put(-27,65){$\Gamma_0(9)$}
\put(-27,55){cusp width:}

\put(47,55){$9$} 
\put(37,65){$t=\infty$}
\put(50,50){\circle*{3}}
\put(50,50){\line(0,-1){30}}
\put(50,20){\circle*{3}}
\put(45,10){$27$}

\put(97,55){$1$}
\put(84,65){$t=-3$}
\put(100,50){\circle*{3}}
\put(100,50){\line(0,-1){30}}
\put(100,20){\circle*{3}}
\put(97,10){$3$}

\put(147,55){$1$}
\put(133,65){$t=-3\omega$}
\put(150,50){\circle*{3}}
\put(150,50){\line(-1,-3){10}}
\put(150,50){\line(0,-1){30}}
\put(150,50){\line(1,-3){10}}
\put(140,20){\circle*{3}}
\put(150,20){\circle*{3}}
\put(160,20){\circle*{3}}
\put(137,10){$1$}
\put(147,10){$1$}
\put(157,10){$1$}

\put(197,55){$1$}
\put(183,65){$t=-3\bar{\omega}$}
\put(200,50){\circle*{3}}
\put(200,50){\line(-1,-3){10}}
\put(200,50){\line(1,-3){10}}
\put(200,50){\line(0,-1){30}}
\put(190,20){\circle*{3}}
\put(200,20){\circle*{3}}
\put(210,20){\circle*{3}}
\put(187,10){$1$}
\put(197,10){$1$}
\put(207,10){$1$}

\put(-27,10){cusp width:}
\put(-27,0){$\Gamma(27;27,3,3)$}
\end{picture}
\end{center}

\subsection{Two index $3$ and two index $4$ subgroups of $\Gamma_0(6)$}\mbox{ }\\
The group $\Gamma_0(6)$ we already discussed in Paragraph~\ref{24-2}.
A corresponding modular ellipic surface is
\[
(x+y)(x+1)(y+1)+txy=0,
\]
with at $t=-8$ an $I_1$-fiber, at $t=1$ an $I_2$-fiber, at $t=0$ an
$I_3$-fiber, and at $t=\infty$ an $I_6$-fiber.
The degree $3$ cyclic covering of $X_0(6)$ ramified over
the cusps of width $1$ and $2$ corresponds to $\Gamma_0(2)\cap\Gamma(3)$.

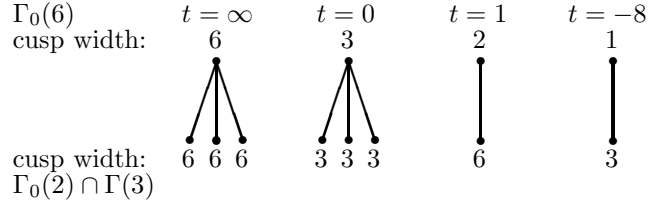
\vspace{\baselineskip}

\begin{center}
\begin{picture}(230,70)(-30,0)
\thicklines
\put(-27,65){$\Gamma_0(6)$}
\put(-27,55){cusp width:}

\put(47,55){$6$} 
\put(37,65){$t=\infty$}
\put(50,50){\circle*{3}}
\put(50,50){\line(0,-1){30}}
\put(50,20){\circle*{3}}
\put(50,50){\line(-1,-3){10}}
\put(50,50){\line(0,-1){30}}
\put(50,50){\line(1,-3){10}}
\put(40,20){\circle*{3}}
\put(50,20){\circle*{3}}
\put(60,20){\circle*{3}}
\put(37,10){$6$}
\put(47,10){$6$}
\put(57,10){$6$}

\put(97,55){$3$}
\put(88,65){$t=0$}
\put(100,50){\circle*{3}}
\put(100,50){\line(-1,-3){10}}
\put(100,50){\line(0,-1){30}}
\put(100,50){\line(1,-3){10}}
\put(90,20){\circle*{3}}
\put(100,20){\circle*{3}}
\put(110,20){\circle*{3}}
\put(87,10){$3$}
\put(97,10){$3$}
\put(107,10){$3$}

\put(147,55){$2$}
\put(138,65){$t=1$}
\put(150,50){\circle*{3}}
\put(150,50){\line(0,-1){30}}
\put(150,20){\circle*{3}}
\put(147,10){$6$}

\put(197,55){$1$}
\put(183,65){$t=-8$}
\put(200,50){\circle*{3}}
\put(200,50){\line(0,-1){30}}
\put(200,20){\circle*{3}}
\put(197,10){$3$}

\put(-27,10){cusp width:}
\put(-27,0){$\Gamma_0(2)\cap \Gamma(3)$}
\end{picture}
\end{center}

\vspace{\baselineskip}

Similarly, the degree $3$ cyclic cover ramified over the cusps
of width $3$ and $6$ yields $\Gamma_0(18)$.

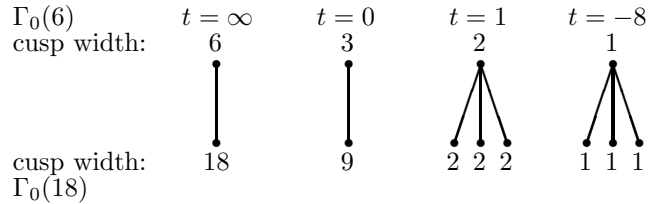
\vspace{\baselineskip}

\begin{center}
\begin{picture}(230,70)(-30,0)
\thicklines
\put(-27,65){$\Gamma_0(6)$}
\put(-27,55){cusp width:}

\put(47,55){$6$} 
\put(37,65){$t=\infty$}
\put(50,50){\circle*{3}}
\put(50,50){\line(0,-1){30}}
\put(50,20){\circle*{3}}
\put(45,10){$18$}

\put(97,55){$3$}
\put(88,65){$t=0$}
\put(100,50){\circle*{3}}
\put(100,50){\line(0,-1){30}}
\put(100,20){\circle*{3}}
\put(97,10){$9$}

\put(147,55){$2$}
\put(138,65){$t=1$}
\put(150,50){\circle*{3}}
\put(150,50){\line(-1,-3){10}}
\put(150,50){\line(0,-1){30}}
\put(150,50){\line(1,-3){10}}
\put(140,20){\circle*{3}}
\put(150,20){\circle*{3}}
\put(160,20){\circle*{3}}
\put(137,10){$2$}
\put(147,10){$2$}
\put(157,10){$2$}

\put(197,55){$1$}
\put(183,65){$t=-8$}
\put(200,50){\circle*{3}}
\put(200,50){\line(-1,-3){10}}
\put(200,50){\line(1,-3){10}}
\put(200,50){\line(0,-1){30}}
\put(190,20){\circle*{3}}
\put(200,20){\circle*{3}}
\put(210,20){\circle*{3}}
\put(187,10){$1$}
\put(197,10){$1$}
\put(207,10){$1$}

\put(-27,10){cusp width:}
\put(-27,0){$\Gamma_0(18)$}
\end{picture}
\end{center}

\vspace{\baselineskip}

Two more subgroups of $\Gamma_0(6)$ appearing in Table~1 have index $4$.
The degree $4$ cyclic covering of $X_0(6)$ ramifying over the cusps
of width $3$ and $1$ corresponds to $\Gamma(12;6,1,2)$.

\vspace{\baselineskip}

\begin{center}
\begin{picture}(230,70)(-30,0)
\thicklines
\put(-27,65){$\Gamma_0(6)$}
\put(-27,55){cusp width:}

\put(47,55){$6$} 
\put(37,65){$t=\infty$}
\put(50,50){\circle*{3}}
\put(50,50){\line(-2,-5){12}}
\put(50,50){\line(-1,-6){5}}
\put(50,50){\line(1,-6){5}}
\put(50,50){\line(2,-5){12}}
\put(38,20){\circle*{3}}
\put(45,20){\circle*{3}}
\put(55,20){\circle*{3}}
\put(62,20){\circle*{3}}
\put(35,10){$6$}
\put(42,10){$6$}
\put(52,10){$6$}
\put(59,10){$6$}

\put(97,55){$3$}
\put(88,65){$t=0$}
\put(100,50){\circle*{3}}
\put(100,50){\line(0,-1){30}}
\put(100,20){\circle*{3}}
\put(95,10){$12$}

\put(147,55){$2$}
\put(138,65){$t=1$}
\put(150,50){\circle*{3}}
\put(150,50){\line(-2,-5){12}}
\put(150,50){\line(-1,-6){5}}
\put(150,50){\line(1,-6){5}}
\put(150,50){\line(2,-5){12}}
\put(138,20){\circle*{3}}
\put(145,20){\circle*{3}}
\put(155,20){\circle*{3}}
\put(162,20){\circle*{3}}
\put(135,10){$2$}
\put(142,10){$2$}
\put(152,10){$2$}
\put(159,10){$2$}

\put(197,55){$1$}
\put(183,65){$t=-8$}
\put(200,50){\circle*{3}}
\put(200,50){\line(0,-1){30}}
\put(200,20){\circle*{3}}
\put(197,10){$4$}

\put(-27,10){cusp width:}
\put(-27,0){$\Gamma(12;6,1,2)$}
\end{picture}
\end{center}

\vspace{\baselineskip}

Finally, the degree $4$ cyclic covering which ramifies over the
cusps of width $6$ and $2$ comes from the group $\Gamma(24;24,2,2)$.
 
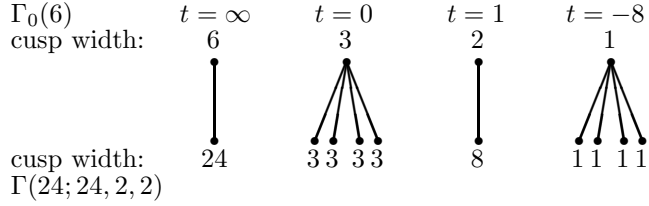
\vspace{\baselineskip}

\begin{center}
\begin{picture}(230,70)(-30,0)
\thicklines
\put(-27,65){$\Gamma_0(6)$}
\put(-27,55){cusp width:}

\put(47,55){$6$} 
\put(37,65){$t=\infty$}
\put(50,50){\circle*{3}}
\put(50,50){\line(0,-1){30}}
\put(50,20){\circle*{3}}
\put(45,10){$24$}

\put(97,55){$3$}
\put(88,65){$t=0$}
\put(100,50){\circle*{3}}
\put(100,50){\line(-2,-5){12}}
\put(100,50){\line(-1,-6){5}}
\put(100,50){\line(1,-6){5}}
\put(100,50){\line(2,-5){12}}
\put(88,20){\circle*{3}}
\put(95,20){\circle*{3}}
\put(105,20){\circle*{3}}
\put(112,20){\circle*{3}}
\put(85,10){$3$}
\put(92,10){$3$}
\put(102,10){$3$}
\put(109,10){$3$}

\put(147,55){$2$}
\put(138,65){$t=1$}
\put(150,50){\circle*{3}}
\put(150,50){\line(0,-1){30}}
\put(150,20){\circle*{3}}
\put(147,10){$8$}

\put(197,55){$1$}
\put(183,65){$t=-8$}
\put(200,50){\circle*{3}}
\put(200,50){\line(-2,-5){12}}
\put(200,50){\line(-1,-6){5}}
\put(200,50){\line(1,-6){5}}
\put(200,50){\line(2,-5){12}}
\put(188,20){\circle*{3}}
\put(195,20){\circle*{3}}
\put(205,20){\circle*{3}}
\put(212,20){\circle*{3}}
\put(185,10){$1$}
\put(192,10){$1$}
\put(202,10){$1$}
\put(209,10){$1$}

\put(-27,10){cusp width:}
\put(-27,0){$\Gamma(24;24,2,2)$}
\end{picture}
\end{center}

\subsection{Two index $4$ subgroups of $\Gamma_0(8)$}\mbox{ }\\
From Paragraph~\ref{24-3} one sees that $\Gamma_0(8)$ gives the
elliptic modular surface
\[
\eta^2=\xi^3+(2-t^2)\xi^2+\xi,
\]
with two $I_1$-fibers at $t=\pm 1$, an $I_2$-fiber at $t=0$ and an
$I_8$-fiber at $t=\infty$.
The groups $\Gamma_1(8)\cap\Gamma(2)$ and $\Gamma(32;32,4,2)$
both have index $4$ in $\Gamma_0(8)$. They correspond to
cyclic coverings of degree $4$ which are ramified over
the cusps of width $1$ for the first group, and over
the cusps of width $2$ and $8$ for the second.
This is summarized in next two diagrams.

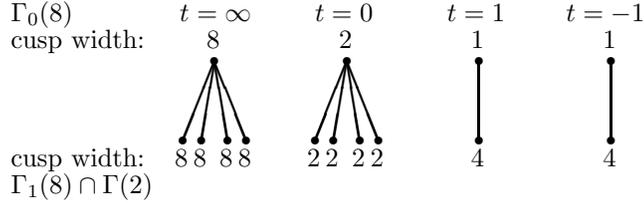
\vspace{\baselineskip}

\begin{center}
\begin{picture}(230,70)(-30,0)
\thicklines
\put(-27,65){$\Gamma_0(8)$}
\put(-27,55){cusp width:}

\put(47,55){$8$} 
\put(37,65){$t=\infty$}
\put(50,50){\circle*{3}}
\put(50,50){\line(-2,-5){12}}
\put(50,50){\line(-1,-6){5}}
\put(50,50){\line(1,-6){5}}
\put(50,50){\line(2,-5){12}}
\put(38,20){\circle*{3}}
\put(45,20){\circle*{3}}
\put(55,20){\circle*{3}}
\put(62,20){\circle*{3}}
\put(35,10){$8$}
\put(42,10){$8$}
\put(52,10){$8$}
\put(59,10){$8$}

\put(97,55){$2$}
\put(88,65){$t=0$}
\put(100,50){\circle*{3}}
\put(100,50){\line(-2,-5){12}}
\put(100,50){\line(-1,-6){5}}
\put(100,50){\line(1,-6){5}}
\put(100,50){\line(2,-5){12}}
\put(88,20){\circle*{3}}
\put(95,20){\circle*{3}}
\put(105,20){\circle*{3}}
\put(112,20){\circle*{3}}
\put(85,10){$2$}
\put(92,10){$2$}
\put(102,10){$2$}
\put(109,10){$2$}

\put(147,55){$1$}
\put(138,65){$t=1$}
\put(150,50){\circle*{3}}
\put(150,50){\line(0,-1){30}}
\put(150,20){\circle*{3}}
\put(147,10){$4$}

\put(197,55){$1$}
\put(183,65){$t=-1$}
\put(200,50){\circle*{3}}
\put(200,50){\line(0,-1){30}}
\put(200,20){\circle*{3}}
\put(197,10){$4$}

\put(-27,10){cusp width:}
\put(-27,0){$\Gamma_1(8)\cap\Gamma(2)$}
\end{picture}
\end{center}

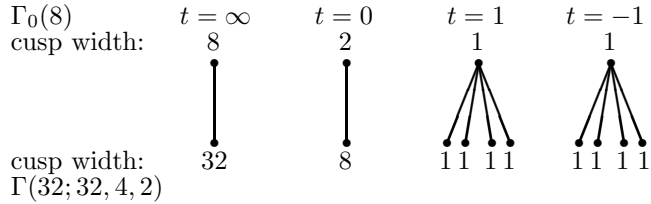
\vspace{\baselineskip}

\begin{center}
\begin{picture}(230,70)(-30,0)
\thicklines
\put(-27,65){$\Gamma_0(8)$}
\put(-27,55){cusp width:}

\put(47,55){$8$} 
\put(37,65){$t=\infty$}
\put(50,50){\circle*{3}}
\put(50,50){\line(0,-1){30}}
\put(50,20){\circle*{3}}
\put(45,10){$32$}

\put(97,55){$2$}
\put(88,65){$t=0$}
\put(100,50){\circle*{3}}
\put(100,50){\line(0,-1){30}}
\put(100,20){\circle*{3}}
\put(97,10){$8$}

\put(147,55){$1$}
\put(138,65){$t=1$}
\put(150,50){\circle*{3}}
\put(150,50){\line(-2,-5){12}}
\put(150,50){\line(-1,-6){5}}
\put(150,50){\line(1,-6){5}}
\put(150,50){\line(2,-5){12}}
\put(138,20){\circle*{3}}
\put(145,20){\circle*{3}}
\put(155,20){\circle*{3}}
\put(162,20){\circle*{3}}
\put(135,10){$1$}
\put(142,10){$1$}
\put(152,10){$1$}
\put(159,10){$1$}

\put(197,55){$1$}
\put(183,65){$t=-1$}
\put(200,50){\circle*{3}}
\put(200,50){\line(-2,-5){12}}
\put(200,50){\line(-1,-6){5}}
\put(200,50){\line(1,-6){5}}
\put(200,50){\line(2,-5){12}}
\put(188,20){\circle*{3}}
\put(195,20){\circle*{3}}
\put(205,20){\circle*{3}}
\put(212,20){\circle*{3}}
\put(185,10){$1$}
\put(192,10){$1$}
\put(202,10){$1$}
\put(210,10){$1$}

\put(-27,10){cusp width:}
\put(-27,0){$\Gamma(32;32,4,2)$}
\end{picture}
\end{center}

\subsection{Two index $4$ subgroups of $\Gamma_0(4)\cap\Gamma(2)$}\mbox{ }\\
As we saw in Paragraph~\ref{24-1}, an equation for the elliptic
modular surface related with $\Gamma_0(4)\cap\Gamma(2)$ is
\
\[ \eta^2=\xi^3+(t^2+1)\xi^2+t^2\xi,
\]
having $I_2$-fibers at $t=\pm 1$ and  $I_4$-fibers at $t=0,\infty$.
The group $\Gamma(8;2,1,2)$ comes from a covering ramified over the
cusps of width $2$, and the group $\Gamma(16;8,2,2)$ comes from a
covering ramified over the cusps of width $4$. We express this as follows.

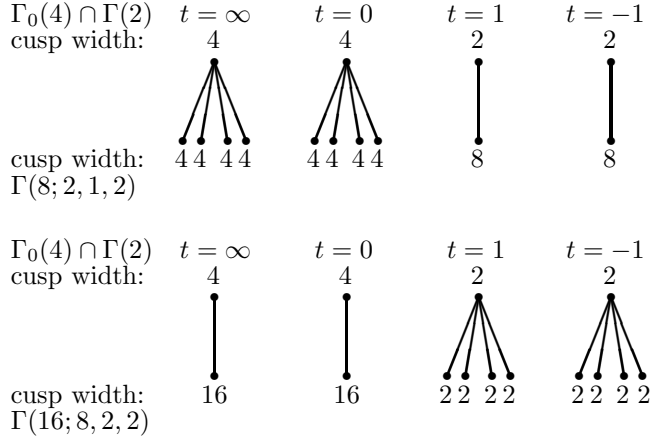
\vspace{\baselineskip}

\begin{center}
\begin{picture}(230,70)(-30,0)
\thicklines
\put(-27,65){$\Gamma_0(4)\cap\Gamma(2)$}
\put(-27,55){cusp width:}

\put(47,55){$4$} 
\put(37,65){$t=\infty$}
\put(50,50){\circle*{3}}
\put(50,50){\line(-2,-5){12}}
\put(50,50){\line(-1,-6){5}}
\put(50,50){\line(1,-6){5}}
\put(50,50){\line(2,-5){12}}
\put(38,20){\circle*{3}}
\put(45,20){\circle*{3}}
\put(55,20){\circle*{3}}
\put(62,20){\circle*{3}}
\put(35,10){$4$}
\put(42,10){$4$}
\put(52,10){$4$}
\put(59,10){$4$}

\put(97,55){$4$}
\put(88,65){$t=0$}
\put(100,50){\circle*{3}}
\put(100,50){\line(-2,-5){12}}
\put(100,50){\line(-1,-6){5}}
\put(100,50){\line(1,-6){5}}
\put(100,50){\line(2,-5){12}}
\put(88,20){\circle*{3}}
\put(95,20){\circle*{3}}
\put(105,20){\circle*{3}}
\put(112,20){\circle*{3}}
\put(85,10){$4$}
\put(92,10){$4$}
\put(102,10){$4$}
\put(109,10){$4$}

\put(147,55){$2$}
\put(138,65){$t=1$}
\put(150,50){\circle*{3}}
\put(150,50){\line(0,-1){30}}
\put(150,20){\circle*{3}}
\put(147,10){$8$}

\put(197,55){$2$}
\put(183,65){$t=-1$}
\put(200,50){\circle*{3}}
\put(200,50){\line(0,-1){30}}
\put(200,20){\circle*{3}}
\put(197,10){$8$}

\put(-27,10){cusp width:}
\put(-27,0){$\Gamma(8;2,1,2)$}
\end{picture}
\end{center}

\vspace{\baselineskip}

\begin{center}
\begin{picture}(230,70)(-30,0)
\thicklines
\put(-27,65){$\Gamma_0(4)\cap\Gamma(2)$}
\put(-27,55){cusp width:}

\put(47,55){$4$} 
\put(37,65){$t=\infty$}
\put(50,50){\circle*{3}}
\put(50,50){\line(0,-1){30}}
\put(50,20){\circle*{3}}
\put(45,10){$16$}

\put(97,55){$4$}
\put(88,65){$t=0$}
\put(100,50){\circle*{3}}
\put(100,50){\line(0,-1){30}}
\put(100,20){\circle*{3}}
\put(95,10){$16$}

\put(147,55){$2$}
\put(138,65){$t=1$}
\put(150,50){\circle*{3}}
\put(150,50){\line(-2,-5){12}}
\put(150,50){\line(-1,-6){5}}
\put(150,50){\line(1,-6){5}}
\put(150,50){\line(2,-5){12}}
\put(138,20){\circle*{3}}
\put(145,20){\circle*{3}}
\put(155,20){\circle*{3}}
\put(162,20){\circle*{3}}
\put(135,10){$2$}
\put(142,10){$2$}
\put(152,10){$2$}
\put(159,10){$2$}

\put(197,55){$2$}
\put(183,65){$t=-1$}
\put(200,50){\circle*{3}}
\put(200,50){\line(-2,-5){12}}
\put(200,50){\line(-1,-6){5}}
\put(200,50){\line(1,-6){5}}
\put(200,50){\line(2,-5){12}}
\put(188,20){\circle*{3}}
\put(195,20){\circle*{3}}
\put(205,20){\circle*{3}}
\put(212,20){\circle*{3}}
\put(185,10){$2$}
\put(192,10){$2$}
\put(202,10){$2$}
\put(210,10){$2$}

\put(-27,10){cusp width:}
\put(-27,0){$\Gamma(16;8,2,2)$}
\end{picture}
\end{center}

\subsection{Two index $5$ subgroups of $\Gamma_1(5)$}\mbox{ }\\
The Beauville elliptic modular surface for $\Gamma_1(5)$ is
given by
\[ x(x-1)(y-1)+ty(x-y)=0. \]
A corresponding Weierstrass equation is obtained using the
coordinate change $\xi:=4t/x$ and $\eta:=4t(2ty-tx-x^2+x)/x^2$;
the resulting equation is
\[
\eta^2=\xi^3+(t^2-6t+1)\xi^2+8t(t-1)\xi+16t^2.
\]
This defines a surface with singular fibers at $t=0,\infty$ (of type $I_5$)
and two fibers of type $I_1$
at $t=\alpha:=\frac{11}{2}+\frac{5}{2}\sqrt{5}$ and at
$t=\beta:=\frac{11}{2}-\frac{5}{2}\sqrt{5}$.

The subgroup $\Gamma(5)$ corresponds to the cyclic quintic covering
of $X_1(5)$ which ramifies over the two cusps of width $1$. 
In terms of the coordinate $t$ which we use on $X_1(5)$,
this covering can be expressed as
\[
s\mapsto t:= \frac{(s+5)(11s^4+70s^3+200s^2+250s+125)}{2s(s^4+50s^2+125)}.
\]
In a diagram, the covering yields the following.

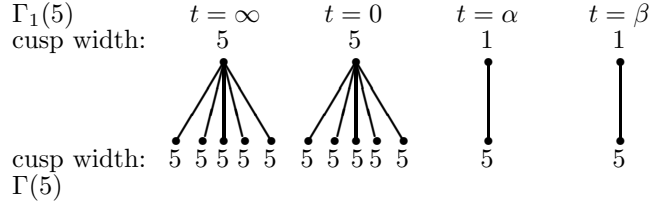
\vspace{\baselineskip}

\begin{center}
\begin{picture}(230,70)(-40,0)
\thicklines
\put(-30,65){$\Gamma_1(5)$}
\put(-30,55){cusp width:}

\put(47,55){$5$} 
\put(37,65){$t=\infty$}
\put(50,50){\circle*{3}}
\put(50,50){\line(0,-1){30}}
\put(50,50){\line(-3,-5){18}}
\put(50,50){\line(-1,-4){7}}
\put(50,50){\line(1,-4){7}}
\put(50,50){\line(3,-5){18}}
\put(32,20){\circle*{3}}
\put(42,20){\circle*{3}}
\put(50,20){\circle*{3}}
\put(58,20){\circle*{3}}
\put(68,20){\circle*{3}}
\put(29,10){$5$}
\put(39,10){$5$}
\put(55,10){$5$}
\put(65,10){$5$}
\put(47,10){$5$}

\put(97,55){$5$}
\put(88,65){$t=0$}
\put(100,50){\circle*{3}}
\put(100,50){\line(0,-1){30}}
\put(100,50){\line(-3,-5){18}}
\put(100,50){\line(-1,-4){7}}
\put(100,50){\line(1,-4){7}}
\put(100,50){\line(3,-5){18}}
\put(100,20){\circle*{3}}
\put(82,20){\circle*{3}}
\put(92,20){\circle*{3}}
\put(108,20){\circle*{3}}
\put(118,20){\circle*{3}}
\put(79,10){$5$}
\put(98,10){$5$}
\put(105,10){$5$}
\put(115,10){$5$}
\put(89,10){$5$}

\put(147,55){$1$}
\put(138,65){$t=\alpha$}
\put(150,50){\circle*{3}}
\put(150,50){\line(0,-1){30}}
\put(150,20){\circle*{3}}
\put(147,10){$5$}

\put(197,55){$1$}
\put(188,65){$t=\beta$}
\put(200,50){\circle*{3}}
\put(200,50){\line(0,-1){30}}
\put(200,20){\circle*{3}}
\put(197,10){$5$}

\put(-30,10){cusp width:}
\put(-30,0){$\Gamma(5)$}
\end{picture}
\end{center}

\medskip
The remaining subgroup $\Gamma_0(25)\cap\Gamma_1(5)$ of $\Gamma_1(5)$ 
corresponds to the cyclic quintic covering
of $X_1(5)$ which ramifies over the two cusps of width $5$. 

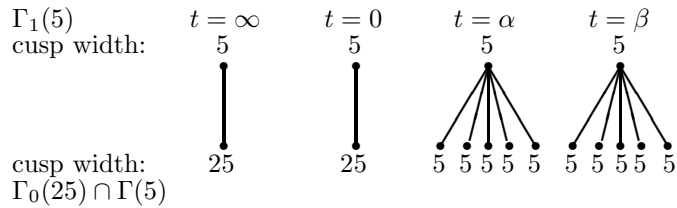
\vspace{\baselineskip}

\begin{center}
\begin{picture}(230,70)(-40,0)
\thicklines
\put(-30,65){$\Gamma_1(5)$}
\put(-30,55){cusp width:}

\put(147,55){$5$} 
\put(137,65){$t=\alpha$}
\put(150,50){\circle*{3}}
\put(150,50){\line(0,-1){30}}
\put(150,50){\line(-3,-5){18}}
\put(150,50){\line(-1,-4){7}}
\put(150,50){\line(1,-4){7}}
\put(150,50){\line(3,-5){18}}
\put(132,20){\circle*{3}}
\put(142,20){\circle*{3}}
\put(150,20){\circle*{3}}
\put(158,20){\circle*{3}}
\put(168,20){\circle*{3}}
\put(129,10){$5$}
\put(139,10){$5$}
\put(155,10){$5$}
\put(165,10){$5$}
\put(147,10){$5$}

\put(197,55){$5$}
\put(188,65){$t=\beta$}
\put(200,50){\circle*{3}}
\put(200,50){\line(0,-1){30}}
\put(200,50){\line(-3,-5){18}}
\put(200,50){\line(-1,-4){7}}
\put(200,50){\line(1,-4){7}}
\put(200,50){\line(3,-5){18}}
\put(200,20){\circle*{3}}
\put(182,20){\circle*{3}}
\put(192,20){\circle*{3}}
\put(208,20){\circle*{3}}
\put(218,20){\circle*{3}}
\put(179,10){$5$}
\put(198,10){$5$}
\put(205,10){$5$}
\put(216,10){$5$}
\put(189,10){$5$}

\put(47,55){$5$}
\put(37,65){$t=\infty$}
\put(50,50){\circle*{3}}
\put(50,50){\line(0,-1){30}}
\put(50,20){\circle*{3}}
\put(44,10){$25$}

\put(97,55){$5$}
\put(88,65){$t=0$}
\put(100,50){\circle*{3}}
\put(100,50){\line(0,-1){30}}
\put(100,20){\circle*{3}}
\put(94,10){$25$}

\put(-30,10){cusp width:}
\put(-30,0){$\Gamma_0(25)\cap\Gamma(5)$}
\end{picture}
\end{center}

\subsection{The remaining groups: $\Gamma_1(12)$ and $\Gamma_0(16)\cap\Gamma_1(8)$}
\mbox{ }\\
With respect to the natural morphism attached to inclusion of groups, the modular curves 
associated with
$\Gamma_1(12)$ and with $\Gamma_0(16)\cap\Gamma_1(8)$ are not cyclic coverings of any
of the modular curves appearing in Beauville's list: such a cyclic morphism would
have degree $4$, and then two of the cusps should have a width divisible by $4$
and the others should come in $4$-tuples having the same width. As can
be read off from Table~1, this is not the case.

To find explicit equations for the corresponding elliptic modular surfaces,
we use the inclusions
\[
\Gamma_1(12)\subset \Gamma_0(12)\subset \Gamma_0(6)
\]
and
\[
\Gamma_0(16)\cap\Gamma_1(8)\subset \Gamma_0(16)\subset \Gamma_0(8)
\]
in which each subgroup has index two in the next one.
An explicit covering morphism can be described as follows.
\[\mbox{ for }\;\Gamma_1(12):\quad s\mapsto t:= 1 - (s^2+3)^2/(s^2+1)^2;
\]
\[
\mbox{ for }\;\Gamma_0(16)\cap\Gamma_1(8):\quad s\mapsto t:= (s^2-1)^2/(s^2+1)^2.
\]
The associated covering data is expressed in the following two diagrams.

\medskip

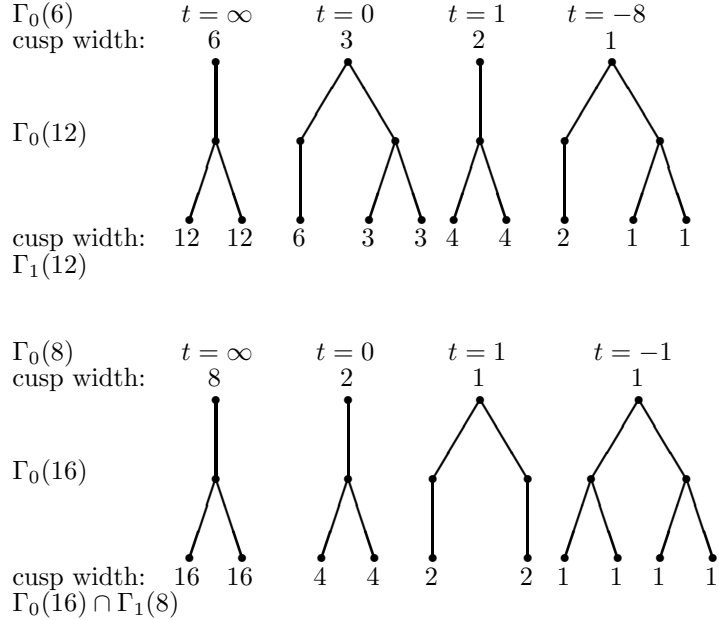
\begin{center}
\begin{picture}(260,110)(-30,-30)
\thicklines
\put(-27,65){$\Gamma_0(6)$}
\put(-27,55){cusp width:}

\put(47,55){$6$} 
\put(37,65){$t=\infty$}
\put(50,50){\circle*{3}}
\put(50,50){\line(0,-1){30}}
\put(50,20){\circle*{3}}
\put(50,20){\line(-1,-3){10}}
\put(50,20){\line(1,-3){10}}
\put(40,-10){\circle*{3}}
\put(60,-10){\circle*{3}}
\put(34,-20){$12$}
\put(54,-20){$12$}

\put(147,55){$2$}
\put(138,65){$t=1$}
\put(150,50){\circle*{3}}
\put(150,50){\line(0,-1){30}}
\put(150,20){\circle*{3}}
\put(150,20){\line(-1,-3){10}}
\put(150,20){\line(1,-3){10}}
\put(140,-10){\circle*{3}}
\put(160,-10){\circle*{3}}
\put(137,-20){$4$}
\put(157,-20){$4$}

\put(97,55){$3$}
\put(88,65){$t=0$}
\put(100,50){\circle*{3}}
\put(100,50){\line(-3,-5){18}}
\put(100,50){\line(3,-5){18}}
\put(82,20){\circle*{3}}
\put(118,20){\circle*{3}}
\put(82,20){\line(0,-1){30}}
\put(82,-10){\circle*{3}}
\put(118,20){\line(-1,-3){10}}
\put(118,20){\line(1,-3){10}}
\put(108,-10){\circle*{3}}
\put(128,-10){\circle*{3}}
\put(79,-20){$6$}
\put(105,-20){$3$}
\put(125,-20){$3$}

\put(197,55){$1$}
\put(183,65){$t=-8$}
\put(200,50){\circle*{3}}
\put(200,50){\line(-3,-5){18}}
\put(200,50){\line(3,-5){18}}
\put(182,20){\circle*{3}}
\put(218,20){\circle*{3}}
\put(182,20){\line(0,-1){30}}
\put(182,-10){\circle*{3}}
\put(218,20){\line(-1,-3){10}}
\put(218,20){\line(1,-3){10}}
\put(208,-10){\circle*{3}}
\put(228,-10){\circle*{3}}
\put(179,-20){$2$}
\put(205,-20){$1$}
\put(225,-20){$1$}

\put(-27,20){$\Gamma_0(12)$}

\put(-27,-20){cusp width:}
\put(-27,-30){$\Gamma_1(12)$}
\end{picture}
\end{center}

\vspace{\baselineskip}

\begin{center}
\begin{picture}(260,110)(-30,-30)
\thicklines
\put(-27,65){$\Gamma_0(8)$}
\put(-27,55){cusp width:}

\put(47,55){$8$} 
\put(37,65){$t=\infty$}
\put(50,50){\circle*{3}}
\put(50,50){\line(0,-1){30}}
\put(50,20){\circle*{3}}
\put(50,20){\line(-1,-3){10}}
\put(50,20){\line(1,-3){10}}
\put(40,-10){\circle*{3}}
\put(60,-10){\circle*{3}}
\put(34,-20){$16$}
\put(54,-20){$16$}

\put(97,55){$2$}
\put(88,65){$t=0$}
\put(100,50){\circle*{3}}
\put(100,50){\line(0,-1){30}}
\put(100,20){\circle*{3}}
\put(100,20){\line(-1,-3){10}}
\put(100,20){\line(1,-3){10}}
\put(90,-10){\circle*{3}}
\put(110,-10){\circle*{3}}
\put(87,-20){$4$}
\put(107,-20){$4$}

\put(147,55){$1$}
\put(138,65){$t=1$}
\put(150,50){\circle*{3}}
\put(150,50){\line(-3,-5){18}}
\put(150,50){\line(3,-5){18}}
\put(132,20){\circle*{3}}
\put(168,20){\circle*{3}}
\put(132,20){\line(0,-1){30}}
\put(132,-10){\circle*{3}}
\put(168,20){\line(0,-1){30}}
\put(168,-10){\circle*{3}}
\put(129,-20){$2$}
\put(165,-20){$2$}

\put(207,55){$1$}
\put(193,65){$t=-1$}
\put(210,50){\circle*{3}}
\put(210,50){\line(-3,-5){18}}
\put(210,50){\line(3,-5){18}}
\put(192,20){\circle*{3}}
\put(228,20){\circle*{3}}
\put(192,20){\line(-1,-3){10}}
\put(192,20){\line(1,-3){10}}
\put(182,-10){\circle*{3}}
\put(202,-10){\circle*{3}}
\put(228,20){\line(-1,-3){10}}
\put(228,20){\line(1,-3){10}}
\put(218,-10){\circle*{3}}
\put(238,-10){\circle*{3}}
\put(179,-20){$1$}
\put(199,-20){$1$}
\put(215,-20){$1$}
\put(235,-20){$1$}

\put(-27,20){$\Gamma_0(16)$}

\put(-27,-20){cusp width:}
\put(-27,-30){$\Gamma_0(16)\cap\Gamma_1(8)$}
\end{picture}
\end{center}

\section{Some remarks}
\label{remarks}

As remarked at the end of Section~\ref{elmod}, the elliptic modular surfaces of
which we found equations are semi-stable and extremal. As is well known (in fact,
the argument in the last paragraph of Section~\ref{elmod} shows this), a
semi-stable elliptic surface over $\PP^1$ of geometric genus $p_g$ is extremal
if and only if the number of singular fibers equals $2p_g+4$.
A result of Nori \cite{No}, recently also proven by Kloosterman \cite{Kl03}
implies for the special case of semi-stable elliptic surfaces that such a
surface is extremal precisely when its $j$-invariant is unramified outside $0, 1728$,
and $\infty$, and all points in $j^{-1}(0)$ have ramification index $3$, and
all point in $j^{-1}(1728)$ have ramification index $2$. 
The degree of this $j$-map will be $12(p_g+1)$, and the types $I_{n_i}$ of the
singular fibers of the elliptic surface can be read off from the fact that
the $n_i$ are the ramification indices of the points in $j^{-1}(\infty)$.

For the case $p_g=0$, i.e., rational elliptic surfaces, Beauville \cite{Be} presented
all semi-stable extremal ones; see also Table~2 above. The next case $p_g=1$,
i.e., elliptic K3 surfaces, Miranda and Persson \cite{MP} showed that there exist
$112$ possible configurations of singular fibers of semi-stable extremal ones.
In the Paragraphs~\ref{24-1}, \ref{24-2}, and \ref{24-3} we showed
that $8$ of them are obtained as quadratic base changes from the
Beauville surfaces. In fact, by choosing different cusps as
ramification points of the quadratic map, one finds equations for
$10$ more such semi-stable extremal elliptic K3 surfaces. Namely, also
the ones with fibre configuration (using the notation from \cite{MP})
\[
[1,1,1,1,2,18],\quad [1,1,1,1,10,10], \quad [1,1,2,2,6,12],\]
\[ [1,1,2,2,9,9],\quad [1,1,2,5,5,10], \quad [1,1,4,6,6,6],\]
\[[2,2,2,3,3,12],\quad [2,2,5,5,5,5],\quad [2,3,3,4,6,6],\quad [3,3,3,3,6,6].\]

If one adds the elliptic modular surface corresponding to $\Gamma_1(7)$ to this
(see Table~3), plus the case $[1,1,1,1,1,19]$ constructed by Iron and Shioda,
then in total $20$ of the $112$ possible cases have been given by explicit
equations. Shioda's method \cite{Sh2} seems unfit for producing more
semi-stable extremal elliptic K3 surfaces. It seems not unlikely, however,
that Iron's method \cite{I} could give several more examples.
The method of Shioda (based on work of Stothers) does however show the
existence of certain semi-stable extremal elliptic surfaces of general type.
More precisely, it shows that for every integer $m>0$ a semi-stabe extremal
elliptic surface exists, with $p_g=m-1$, having $2m+1$ fibers of type $I_1$
and $1$ fiber of type $I_{10m-1}$. 

Actually, it should be possible to obtain a few more equations by using
base changes from easier (rational) surfaces. To illustrate this point,
start with a rational elliptic surface with two fibers of type $II$,
one of type $I_1$, and one of type $I_7$. The existence of this follows from
the existence of the corresponding $j$-invariant: it suffices to show
that a rational function $j(t)$ of the form
\[
j(t):=\lambda \frac{(t-a)^3(t^2+bt+c)}{t^7(t-d)}
\]
can be constructed, with $a,d$ distinct and different from $0$,
and $t^2+bt+c$ having two different zeroes both different from
any of $a,d,0$, and $j(t)=1728$ having only roots with even multiplicity.
Then an elliptic surface with this $j$-invariant exists, having fibers
of type $II$ over the roots of $t^2+bt+c=0$, and an $I_7$-fiber over
$t=0$ and an $I_1$ fiber over $t=d$ and no other singular fibers.
Now taking a cyclic base change of degree $3$ ramified over
the zeroes of $t^2+bt+c$, and then twisting the resulting surface over
the quadratic extension ramified only at these two zeroes,
results in a semi-stable surface with $[1,1,1,7,7,7]$-configuration.

Nevertheless, for most cases such easy cyclic base changes will not
exist.
Therefore, we consider the problem of finding explicit equations for the
remaining $112-20=92$ cases in the table of Miranda and Persson wide open.

\end{document}